\documentclass[11pt,twoside]{amsart}
\usepackage{amsmath, amsthm, amscd, amsfonts, amssymb, graphicx, color}
\usepackage[bookmarksnumbered, colorlinks, plainpages]{hyperref}
\usepackage{pgf,tikz}
\usetikzlibrary{arrows}
\usepackage[utf8]{inputenc}
\usepackage{pstricks-add}
\input{mathrsfs.sty}
\newtheorem{thm}{Theorem}[section]
\newtheorem{cor}[thm]{Corollary}
\newtheorem{lem}[thm]{Lemma}
\newtheorem{prop}[thm]{Proposition}
\newtheorem{defn}[thm]{Definition}

\newtheorem{example}[thm]{\bf{Example}}

\numberwithin{equation}{section}

\begin{document}


\title{Unmixed $d$-uniform $r$-partite Hypergraphs}
\author{Reza jafarpour-Golzari and Rashid Zaare-Nahandi}
\address{Department of Mathematics, Institute for Advanced Studies
in Basic Science (IASBS), P.O.Box 45195-1159, Zanjan, Iran}

\email{r.golzary@iasbs.ac.ir}

\address{Department of Mathematics, Institute for Advanced Studies
in Basic Science (IASBS), P.O.Box 45195-1159, Zanjan, Iran}

\email{rashidzn@iasbs.ac.ir}

\thanks{{\scriptsize
\hskip -0.4 true cm MSC(2010): Primary: 5E40; Secondary: 5C65.
\\
\newline\indent{\scriptsize}}}

\maketitle


\begin{abstract}
In this paper, we characterize all unmixed $d$-uniform $r$-partite hypergraphs under a certain condition. Also we give a necessary condition for unmixedness in $d$-uniform hypergraphs with a perfect matching of size $n$. Finally we give a sufficient condition for unmixednes in $d$-uniform hypergraphs with a perfect matching. 
\end{abstract}

\vskip 0.2 true cm


\section{\bf Introduction}
\vskip 0.4 true cm

Unmixedness is one of the most important concepts in theory of graphs and hypergraphs with nice and interesting algebraic and geometric interpretations (for instance see \cite{Her 2}, \cite{Kia}, \cite{Mor}, \cite{Sta}, \cite{Vil 1}, \cite{Zar}). According to this, characterization of special classes of unmixed graphs has been noteworthy in recent years. G. Ravindra in \cite{Rav} and R. H. Villarreal in \cite{Vil 2} have characterized all unmixed bipartite graphs  independently. H. Haghighi in \cite{Hag} has given a characterization for unmixed tripartite graphs under a certain condition, and recently R. Jafarpour-Golzari and R. Zaare-Nahandi in \cite{Jaf} have generalized Haghighi's result for unmixed $r$-partite graphs. On characterization of unmixed $r$-partite hypergraphs, almost no study has been done. Only in \cite{Kia} a classification of a very special class of unmixed multipartite hypergraphs has been provided.
  
In this paper we give a characterization of all unmixed $d$-uniform $r$-partite hypergraphs under a certain condition which we name it ($\ast\ast$). Also we give necessary or sufficient conditions for unmixedness in more general classes of $d$-uniform hypergraphs (Propositions 3.6, 3.8).
  

\section{\bf {\bf \em{\bf Preliminaries}}}
\vskip 0.4 true cm

In the sequel, we use \cite{Wes} and \cite{Ber} for terminology and notations on graphs and hypergraphs respectively.

Let $G=(V, E)$ be a simple finite graph. For $x, y\in V$, $x\sim y$ means that $x$ and $y$ are adjacent. A subset $M$ of $V$ is said to be independent if for every $x, y\in M$, $x\nsim y$. A vertex cover for $G$ is a subset $C$ of $V$ such that every edge of $G$, intersects $C$. A vertex cover $C$ is minimal whenever there is no any pure subset of $C$ which is a vertex cover. $G$ is called is called unmixed if all minimal vertex cover of $G$ have the same number of elements. A subset $Q$ of $V$ is said to be a clique if for every two distinct vertices $x, y\in Q$, $x\thicksim y$.  

A hypergraph $\mathcal{H}$ on a finite nonempty set $V$ is a set of nonempty subsets of $V$ such that $\bigcup_{e\in\mathcal{H}}  e=V$. The elements of $V$ are called vertices and each element of $\mathcal{H}$ is said a hyperedge. We denote by $V(\mathcal{H})$ and $E(\mathcal{H})$, the sets of vertices and hyperedges of $\mathcal{H}$ respectively. A hypergraph is said to be simple hypergraph or clutter if non of its two distinct hyperedges contains another. The hypergraph $\mathcal{H}$ is called $d$-uniform (or $d$-graph), if all its hyperedges have the same cardinality $d$.

\begin{defn}
 An $r$-partite ($r\geq 2$) hypergraph $\mathcal{H}$, is a hypergraph which $V(\mathcal{H})$ can be partitioned to $r$ subsets such that for every two vertices $x, y$ in one part, $x, y$ do not lie in any hyperedge. Such a partition of $V(\mathcal{H})$ is called an $r$-partition of $\mathcal{H}$. If $r=2, 3$, the $r$-partite hypergraph is said to be bipartite and tripartite respectivily. 
\end{defn}

In the hypergraph $\mathcal{H}$, two vertices $x, y$ are said to be adjacent if there is a hyperedge containing $x$ and $y$. We say that a hyperedge $e$ is adjacent with a vertex $x$ if $x\in e$. For a vertex $x$ of $\mathcal{H}$, the neighborhood of $x$, denoted by $N(x)$, is the set of all vertices which are adjacent to $x$.

A subset $M$ of $V(\mathcal{H})$ is called independent if it dose not contain any hyperedge. An independent set $M$ of $\mathcal{H}$ is said to be maximal whenever it is not strictly contained in any other independent set. A subset $C$ of $V(\mathcal{H})$ is called a vertex cover, if every hyperedge of $\mathcal{H}$ intersects it. A vertex cover is said to be minimal if there is no any pure subset of it which is also a vertex cover. It is clear that every maximal independent set is complement of a minimal vertex cover and vice versa.

\begin{defn}
The hypergraph $\mathcal{H}$ is said to be unmixed if all minimal vertex covers of $\mathcal{H}$ have the same cardinality.
\end{defn}

A matching in a hypergraph $\mathcal{H}$ is a set of hyperedges which are disjoint pairwise. A perfect matching is a matching such that every vertex of $\mathcal{H}$ lies in at least one of its elements.

Let $\{1, \ldots , n\}$ is denoted by $[n]$. A simplicial complex on $[n]$ is a set $\Delta$ of subsets of $[n]$ such that (a) $\{x\}\in\Delta$, for every $x\in[n]$, (b) if $F\in\Delta$ and $G\subseteq F$, then $G\in\Delta$. Each element of $\Delta$ is said to be a face. Dimension of the face $F$, denoted by dim $F$, is defined as $|F|-1$ and dimension of $\Delta$ is the maximum of dimentions of its faces.

Let $\Delta$ be a simplicial complex on $[n]$ and $S$ be a nonempty set of subsets of $[n]$ such that $\bigcup_{N\in{S}} N=[n]$. The simplicial complex generated by $S$ is the set of all subsets of elements of $S$.

For a simplicial complex $\Delta$ or a hypergraph $\mathcal{H}$ on $[n]$, and for $r\geq 0$, the $r$-skeleton of $\Delta$ or $\mathcal{H}$, is the set of all faces of $\Delta$ whose dimension is at most $r$ or all subsets of hyperedges of $\mathcal{H}$ with cardinality not exceeding $r+1$, respectively.

\begin{defn}
Let $\mathcal{H}$ be a $d$-uniform ($d\geq 2$) hypergraph. A $(d-1)$-subset of a hyperedge is called a submaximal edge, and the set of all submaximal edges is denoted by $SE(\mathcal{H})$.
\end{defn}

For $\mathfrak{e}\in SE(\mathcal{H})$, the set $\{v\in V(\mathcal{H})\mid \mathfrak{e}\cup\{v\}\in E(\mathcal{H})\}$, is denoted by $N(\mathfrak{e})$. If $v\in N(\mathfrak{e})$, we write $\mathfrak{e}\sim v$.

In a $d$-uniform hypergraph $\mathcal{H}$, a clique is a subset $W$ of $V(\mathcal{H})$ such that every its subset of size $d$, is a hyperedge in $\mathcal{H}$. 
  
 
 \section{\bf {\bf \em{\bf Unmixed hypergraphs}}}
 \vskip 0.4 true cm
 
 Let $\mathcal{H}$ is a $d$-uniform $r$-partite hypergraph with $2\leq d\leq r$. We say that $\mathcal{H}$ satisfies the condition ($\ast\ast$) for $r\geq 2$, if $\mathcal{H}$ can be partitioned to $r$ parts $V_{i}=\{x_{1i}, \cdots , x_{ni}\}$, $1\leq i\leq r$, such that $\{x_{j1}, x_{j2} \ldots , x_{jr}\}$ is a clique for every $1\leq j\leq n$.
 
 The authors in \cite{Jaf} have presented a necessary and sufficient condition for unmixedness of an $r$-partite graph which satisfies the following condition $(\ast)$ for $r\geq 2$.\\
 We say a graph $G$ satisfies the condition $(\ast )$ for an integer $r\geq 2$, if $G$ can be partitioned to $r$ parts
  $V_{i}=\{x_{1i},
  \ldots , x_{ni}\}$, $1\leq i\leq r$, such that for all $1\leq j\leq n$, $\{x_{j1}, \ldots , x_{jr}\}$ is a clique.
  
Let $\mathcal{H}$ be a $d$-uniform $r$-partite hypergraph ($2\leq d\leq r$) on $[n]$ which satisfies the condition ($\ast\ast$) for $r\geq 2$. Then the 1-skeleton of $\mathcal{H}$ is an r-partite graph which satisfies the condition ($\ast)$ for $r$. But in general, the unmixedness of a hypergraph and its
1-skeleton, are two independent facts, as the following example exhibits. 

\begin{example}\rm
The following clutter is not unmixed while its 1-skeleton is unmixed as a graph.

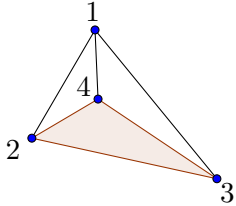
\definecolor{qqqqff}{rgb}{0.,0.,1.}
\begin{center}
\definecolor{zzttqq}{rgb}{0.6,0.2,0.}
\definecolor{qqqqff}{rgb}{0.,0.,1.}
\begin{tikzpicture}[line cap=round,line join=round,>=triangle 45,x=1.0cm,y=1.0cm]
\clip(-0.8,3.5) rectangle (3.34,6.6);
\fill[color=zzttqq,fill=zzttqq,fill opacity=0.1] (0.66,5.08) -- (-0.22,4.56) -- (2.24,4.02) -- cycle;
\draw [color=zzttqq] (0.66,5.08)-- (-0.22,4.56);
\draw [color=zzttqq] (-0.22,4.56)-- (2.24,4.02);
\draw [color=zzttqq] (2.24,4.02)-- (0.66,5.08);
\draw (-0.22,4.56)-- (0.62,6.);
\draw (0.62,6.)-- (2.24,4.02);
\draw (0.66,5.08)-- (0.62,6.);
\draw (0.36,6.5) node[anchor=north west] {1};
\draw (-0.7,4.66) node[anchor=north west] {2};
\draw (2.15,4.10) node[anchor=north west] {3};
\draw (0.24,5.5) node[anchor=north west] {4};
\begin{scriptsize}
\draw [fill=qqqqff] (0.66,5.08) circle (1.5pt);
\draw [fill=qqqqff] (-0.22,4.56) circle (1.5pt);
\draw [fill=qqqqff] (2.24,4.02) circle (1.5pt);
\draw [fill=qqqqff] (0.62,6.) circle (1.5pt);
\end{scriptsize}
\end{tikzpicture}

 \end{center}
 \label{graph 1}

Note that the sets $\{2, 4, 3\}$ and $\{1, 4\}$ are two minimal vertex covers with different sizes for the hypergraph.

Conversely, the following clutter is unmixed but its 1-skeleton is not.

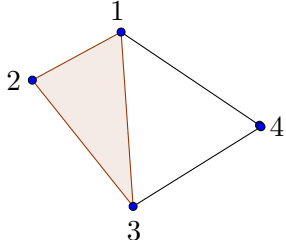
\definecolor{qqqqff}{rgb}{0.,0.,1.}
\begin{center}
\definecolor{zzttqq}{rgb}{0.6,0.2,0.}
\definecolor{qqqqff}{rgb}{0.,0.,1.}
\begin{tikzpicture}[line cap=round,line join=round,>=triangle 45,x=1.0cm,y=1.0cm]
\clip(1.,1.9975) rectangle (5.8,6.2175);
\fill[color=zzttqq,fill=zzttqq,fill opacity=0.1] (3.06,4.76) -- (1.88,4.12) -- (3.22,2.44) -- cycle;
\draw [color=zzttqq] (3.06,4.76)-- (1.88,4.12);
\draw [color=zzttqq] (1.88,4.12)-- (3.22,2.44);
\draw [color=zzttqq] (3.22,2.44)-- (3.06,4.76);
\draw (3.06,4.76)-- (4.9,3.52);
\draw (3.22,2.44)-- (4.92,3.5);
\draw (2.78,5.2975) node[anchor=north west] {1};
\draw (1.4,4.3775) node[anchor=north west] {2};
\draw (3.,2.3775) node[anchor=north west] {3};
\draw (4.9,3.7575) node[anchor=north west] {4};
\begin{scriptsize}
\draw [fill=qqqqff] (3.06,4.76) circle (1.5pt);
\draw [fill=qqqqff] (1.88,4.12) circle (1.5pt);
\draw [fill=qqqqff] (3.22,2.44) circle (1.5pt);
\draw [fill=qqqqff] (4.9,3.52) circle (1.5pt);
\draw [fill=qqqqff] (4.92,3.5) circle (1.5pt);
\end{scriptsize}
\end{tikzpicture}

  \end{center}
 \label{graph 2}

\end{example}

Note that the sets $\{1, 3\}$ and $\{2, 3, 4\}$ are two minimal vertex covers with different sizes for the 1-skeleton.
 
This gives us a motivation for finding a necessary and sufficient condition under which a $d$-uniform $r$-partite ($2\leq d\leq r$) hypergraph satisfying the condition ($\ast\ast$) for $r$, is unmixed.

\

First we prove a lemma.

\begin{lem}
Let $\mathcal{H}$ be a $d$-uniform $r$-partite ($2\leq d\leq r$) hypergraph which satisfies the condition ($\ast\ast$) for $r$. If $\mathcal{H}$ is unmixed, then every minimal vertex cover of $\mathcal{H}$ contains exactly $r-d+1$ elements of each clique $\{x_{j1}, x_{j2} \ldots , x_{jr}\}$.
\begin{proof}
Let $C$ be a minimal vertex cover of $\mathcal{H}$. For every $1\leq q\leq n$, $C$ contains at least $r-d+1$ vertices of the clique $\{x_{q1}, x_{q2} \ldots , x_{qr}\}$, because if $C$ contains at most $r-d$ vertices of that clique, it dose not cover hyperedges on remaining vertices. Therefore a vertex cover must contain at least $n(r-d+1)$ vertices. On the other hand, the set $\bigcup_{i=1}^{r-d+1}V_{i}$ is a minimal vertex cover of $\mathcal{H}$ with $n(r-d+1)$ vertices. This completes the proof.   
\end{proof}
\end{lem}

Now we present the main theorem of this paper.

\begin{thm}
Let $\mathcal{H}$ be a $d$-uniform $r$-partite ($2\leq d\leq r$) hypergraph which satisfies the condition ($\ast\ast$) for $r$. Then $\mathcal{H}$ is unmixed if and only if the following condition holds.\\
For every $1\leq q\leq n$, if $\mathfrak{e}_{1}, \mathfrak{e}_{2} \ldots , \mathfrak{e}_{r-d+2}$ are submaximal edges such that
\[\mathfrak{e}_{1}\thicksim x_{qi_{1}}, \mathfrak{e}_{2}\thicksim x_{qi_{2}}, \ldots , \mathfrak{e}_{r-d+2}\thicksim x_{qi_{r-d+2}}\]
where $i_{1}, i_{2} \ldots , i_{r-d+2}$ are distinct, then the set 
\[\mathfrak{e}_{1}\cup \mathfrak{e}_{2}\cup\ldots \cup \mathfrak{e}_{r-d+2}\]
is not independent.
\begin{proof} 
Let $\mathcal{H}$ be unmixed. We show that the mentioned condition holds. Suppose in contrary
\[\mathfrak{e}_{1}\thicksim x_{qi_{1}}, \mathfrak{e}_{2}\thicksim x_{qi_{2}}, \ldots , \mathfrak{e}_{r-d+2}\thicksim x_{qi_{r-d+2}}\]
where $i_{1}, i_{2} \ldots , i_{r-d+2}$ are distinct but the set
\[F=\mathfrak{e}_{1}\cup \mathfrak{e}_{2}\cup\ldots \cup \mathfrak{e}_{r-d+2}\] 
is independent. Therefore there is a maximal independent set $M$ containing $F$. Since $M$ is a maximal independent set, $C:=V(\mathcal{H})\backslash M$ is a minimal vertex cover of $\mathcal{H}$ which contains no any element of $F$. Since $C$ is a vertex cover of $\mathcal{H}$, $C$ contains the vertices $x_{qi_{1}}, x_{qi_{2}}, \ldots ,x_{qi_{r-d+2}}$. But by Lemma 3.2, $C$ contains exactly $r-d+1$ vertices of every clique, a contradiction.

Conversely, let the mentioned condition holds. We show that $\mathcal{H}$ is unmixed. It is enough to show that every minimal vertex cover of $\mathcal{H}$ contains exactly $r-d+1$ vertices of each clique $\{x_{q1}, x_{q2} \ldots , x_{qr}\}$. Let $C$ be an arbitrary minimal vertex cover and $1\leq q\leq n$. Then $C$ intersects the set $\{x_{q1}, x_{q2} \ldots , x_{qr}\}$ in at least $r-d+1$ elements. Let $C$ intersects the mentioned clique in at least $r-d+2$ elements. Without loss of generality, we assume that this elements are $x_{q1}, x_{q2} \ldots , x_{q(r-d+2)}$. For each $i$, $1\leq i\leq r-d+2$, $x_{qi}$ is in the minimal vertex cover $C$. Then there is a hyperedge $e_{i}$ covered only by $x_{qi}$. That is, $e_{i}\cap C=\{x_{qi}\}$. Suppose that $\mathfrak{e}_{i}=e_{i}\backslash\{x_{qi}\}$. Then the sets $\mathfrak{e}_{i}$ are submaximal edges such that  
\[\mathfrak{e}_{1}\thicksim x_{qi_{1}}, \mathfrak{e}_{2}\thicksim x_{qi_{2}}, \ldots , \mathfrak{e}_{r-d+2}\thicksim x_{qi_{r-d+2}}\]
and $\mathfrak{e}_{1}\cup \mathfrak{e}_{2}\cup\ldots \cup \mathfrak{e}_{r-d+2}$ dose not intersects $C$. But by hypothesis 
\[\mathfrak{e}_{1}\cup \mathfrak{e}_{2}\cup\ldots \cup \mathfrak{e}_{r-d+2}\]
is not independent. That is, it contains a hypergraph $e$ which 
is not covered by $C$, a contradiction.
\end{proof}
\end{thm}

The following theorem of Villarreal on unmixedness of bipartite graphs can be concluded from the Theorem 3.3, where $r=2, d=2$.

\begin{cor}
\cite[Theorem 1.1]{Vil 2} Let $G$ be a bipartite graph without isolated vertices. Then $G$ is unmixed if and only if there is a
bipartition $V_{1}=\{x_{1}, \ldots , x_{g}\}, V_{2}=\{y_{1},
\ldots , y_{g}\}$ of $G$ such that: (a) $\{x_{i}, y_{i}\}\in
E(G)$, for all i, and (b) if $\{x_{i}, y_{j}\}$ and $\{x_{j},
y_{k}\}$ are in $E(G)$, and $i, j, k$ are distinct, then
$\{x_{i}, y_{k}\}\in E(G)$.
\end{cor}

The following theorem can be concluded from theorem 3.3, where d=2.

\begin{cor}
\cite[Theorem 2.3]{Jaf} Let $G$ be an $r$-partite graph which satisfies the
condition $(\ast )$ for $r$. Then $G$ is unmixed if and only if the
following condition hold:
\\For every $1\leq q\leq n$, if there is a set $\{x_{k_{1}s_{1}}, \ldots , x_{k_{r}s_{r}}\}$ such that 
\[x_{k_{1}s_{1}}\thicksim x_{q1}, \ldots , x_{k_{r}s_{r}}\thicksim x_{qr},\]
then the set $\{x_{k_{1}s_{1}}, \ldots , x_{k_{r}s_{r}}\}$ is not independent.
\end{cor}

Now we prove two propositions about $d$-uniform $(d\geq 2)$ hypergraphs by methods used in the proof of Theorem 3.3.

\begin{prop}
Let $\mathcal{H}$ be a d-uniform $(d\geq 2)$ hypergraph on vertex set $\{x_{ji}| \ 1\leq j\leq n, 1\leq i\leq d\}$ with perfect matching  
\[\{\{x_{j1}, x_{j2}, \ldots , x_{jd}\}| \ 1\leq j\leq n\}.\]
If $\mathcal{H}$ is unmixed and has a minimal vertex cover of size $n$, then for every $1\leq q\leq n$, if $\mathfrak{e}_{1}, \mathfrak{e}_{2}$ be two submaximal edges such that
\[\mathfrak{e}_{1}\sim x_{qi_{1}}, \mathfrak{e}_{2}\sim x_{qi_{2}}\]
where $i_{1}, i_{2}$ are distinct, then the set $\mathfrak{e}_{1}\cup \mathfrak{e}_{2}$ is not independent.
\begin{proof}
Let $1\leq q\leq n$ be arbitrary. Let in contrary the set $\mathfrak{e}_{1}\cup \mathfrak{e}_{2}$ is independent . Therefore $\mathfrak{e}_{1}\cup \mathfrak{e}_{2}$ is contained in a maximal independent set $M$. Set $T=V(G)\setminus M$. $T$ is a minimal vertex cover and since it dose not contain any element of $\mathfrak{e}_{1}\cup \mathfrak{e}_{2}$, then $T$ contains $x_{qi_{1}},
x_{qi_{2}}$ and then $T$ is at least of size $n+1$, a contradiction.
\end{proof}
\end{prop}

\begin{example}\rm
In 3-uniform hypergraph
\[\mathcal{H}=\{\{a, b, c\}, \{d, e, f\}, \{g, h, i\}, \{b, g, e\}, \{c, f, h\}\},\]
we have the perfect matching $\{\{a, b, c\}, \{d, e, f\}, \{g, h, i\}\}$ of size 3.

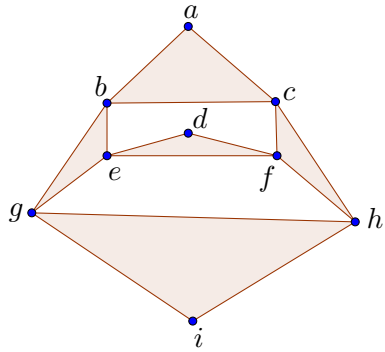
\definecolor{zzttqq}{rgb}{0.6,0.2,0.}
\definecolor{qqqqff}{rgb}{0.,0.,1.}
\begin{center}
\begin{tikzpicture}[line cap=round,line join=round,>=triangle 45,x=1.0cm,y=1.0cm]
\clip(2.,0.) rectangle (8.2,5.7);
\fill[color=zzttqq,fill=zzttqq,fill opacity=0.1] (5.2,4.8) -- (4.12,3.78) -- (6.36,3.8) -- cycle;
\fill[color=zzttqq,fill=zzttqq,fill opacity=0.1] (6.36,3.8) -- (6.38,3.08) -- (7.42,2.2) -- cycle;
\fill[color=zzttqq,fill=zzttqq,fill opacity=0.1] (6.38,3.08) -- (5.2,3.38) -- (4.12,3.08) -- cycle;
\fill[color=zzttqq,fill=zzttqq,fill opacity=0.1] (4.12,3.78) -- (4.12,3.08) -- (3.12,2.32) -- cycle;
\fill[color=zzttqq,fill=zzttqq,fill opacity=0.1] (3.12,2.32) -- (7.42,2.2) -- (5.26,0.88) -- cycle;
\draw [color=zzttqq] (5.2,4.8)-- (4.12,3.78);
\draw [color=zzttqq] (4.12,3.78)-- (6.36,3.8);
\draw [color=zzttqq] (6.36,3.8)-- (5.2,4.8);
\draw [color=zzttqq] (6.36,3.8)-- (6.38,3.08);
\draw [color=zzttqq] (6.38,3.08)-- (7.42,2.2);
\draw [color=zzttqq] (7.42,2.2)-- (6.36,3.8);
\draw [color=zzttqq] (6.38,3.08)-- (5.2,3.38);
\draw [color=zzttqq] (5.2,3.38)-- (4.12,3.08);
\draw [color=zzttqq] (4.12,3.08)-- (6.38,3.08);
\draw (5.14,0.92) node[anchor=north west] {$i$};
\draw (2.68,2.58) node[anchor=north west] {$g$};
\draw [color=zzttqq] (4.12,3.78)-- (4.12,3.08);
\draw [color=zzttqq] (4.12,3.08)-- (3.12,2.32);
\draw [color=zzttqq] (3.12,2.32)-- (4.12,3.78);
\draw [color=zzttqq] (3.12,2.32)-- (7.42,2.2);
\draw [color=zzttqq] (7.42,2.2)-- (5.26,0.88);
\draw [color=zzttqq] (5.26,0.88)-- (3.12,2.32);
\draw (5.12,3.85) node[anchor=north west] {$d$};
\draw (7.44,2.52) node[anchor=north west] {$h$};
\draw (4.,3.06) node[anchor=north west] {$e$};
\draw (6.,3.08) node[anchor=north west] {$f$};
\draw (6.32,4.12) node[anchor=north west] {$c$};
\draw (3.82,4.28) node[anchor=north west] {$b$};
\draw (5.,5.21) node[anchor=north west] {$a$};
\begin{scriptsize}
\draw [fill=qqqqff] (4.12,3.78) circle (1.5pt);
\draw [fill=qqqqff] (5.2,4.8) circle (1.5pt);
\draw [fill=qqqqff] (4.12,3.08) circle (1.5pt);
\draw [fill=qqqqff] (6.36,3.8) circle (1.5pt);
\draw [fill=qqqqff] (6.38,3.08) circle (1.5pt);
\draw [fill=qqqqff] (7.42,2.2) circle (1.5pt);
\draw [fill=qqqqff] (5.2,3.38) circle (1.5pt);
\draw [fill=qqqqff] (5.26,0.88) circle (1.5pt);
\draw [fill=qqqqff] (3.12,2.32) circle (1.5pt);
\end{scriptsize}
\end{tikzpicture}
\end{center}

We show by proposition 3.7 that $\mathcal{H}$ is not unmixed.

Let $\mathcal{H}$ be unmixed (by contrary). $\mathcal{H}$ has the minimal vertex cover $\{b, e, h\}$ of size 3. Now we have 2 hyperedges $\{a, b, c\}$ and $\{b, g, e\}$ in relevence with the hyperedge $\{a, b, c\}$ of perfect matching, but $\{g, e, h, f\}$
is independent, a contradiction.
\end{example}

\begin{prop}
Let $\mathcal{H}$ is a d-uniform $(d\geq 2)$ hypergraph on the vertex set $\{x_{ji}| \ 1\leq j\leq n, 1\leq i\leq d\}$with perfect matching  
\[\{\{x_{j1}, x_{j2}, \ldots , x_{jd}\}| \ 1\leq j\leq n\}.\]
Then a sufficient condition for unmixedness of $\mathcal{H}$ is that for every $1\leq q\leq n$, if $\mathfrak{e}_{1}, \mathfrak{e}_{2}$ be two submaximal edges such that
\[\mathfrak{e}_{1}\sim x_{qi_{1}}, \mathfrak{e}_{2}\sim x_{qi_{2}}\]
where $i_{1}, i_{2}$ are distinct, then $\mathfrak{e}_{1}\cup \mathfrak{e}_{2}$ is not independent.   
\begin{proof}
Let $\mathcal{H}$ satisfies the above condition. We show that $\mathcal{H}$ is unmixed. It is enough to show that every minimal vertex cover of $\mathcal{H}$ contains exactly one element of each hyperedge of the perfect matching. Let $T$ be a minimal vertex cover. $T$ contains at least one element of each hyperedge of perfect matching. Suppose in contrary that T chooses at least two elements from hyperedge $\{x_{q1}, x_{q2} \ldots , x_{qd}\}$. Without loss of generality, let the elements $x_{q1}$ and $x_{q2}$ are chosen. Since $x_{q1}$ is in minimal vertex cover, then there exist at least $d-1$ distinct vertices in $N(x_{q1})$, such that they dose not belong to $T$ and form a hyperedge together with $x_{q1}$. Name the set of these vertices $\mathfrak{e}_{1}$. We have  
\[\mathfrak{e}_{1}\sim x_{q1}.\]
Similarly, Since $x_{q2}$ is in minimal vertex cover, with a similar argument, there is a submaximal edge $\mathfrak{e}_{2}$ consisting of $d-1$ distinct vertices, no one belonging to $T$, such that
\[\mathfrak{e}_{2}\sim x_{q2}.\]

Now $\mathfrak{e}_{1}\cup \mathfrak{e}_{2}$ dose not intersect $T$. But according to hypothesis   
$\mathfrak{e}_{1}\cup \mathfrak{e}_{2}$ is not independent. That is, it contains a hyperedge $e$ which is not covered by $T$, a contradiction. 
\end{proof}
\end{prop}

 
\section{\bf {\bf \em{\bf Edge ideal of unmixed hypergraphs}}}
\vskip 0.4 true cm

In this section, we provide an algebraic interpretation for Theorem 3.3.

\begin{defn}
Let $\mathcal{H}$ be a hypergraph with $V(\mathcal{H})=\{x_{1}, \ldots , x_{m}\}$. Let $K[x_{1}, \ldots , x_{m}]$ be the polynomial ring with indeterminates $x_{1}, \ldots , x_{m}$ and coefficients in a field $K$. For a subset $D=\{x_{i_{1}}, \ldots , x_{i_{r}}\}\subseteq V(\mathcal{H})$, let $X_{D}=x_{i_{1}} \ldots  x_{i_{r}}$. We define the edge ideal of $\mathcal{H}$ to be
\[I(\mathcal{H}):=(X_{e}| \ e\in E(\mathcal{H})).\]  
The quotient ring $K[\mathcal{H}]:=\dfrac{K[x_{1}, \ldots , x_{m}]}{I(\mathcal{H})}$ is called the edge ring of $\mathcal{H}$.
\end{defn}

Let $R$ be a commutative ring. An element $a\in R$ is called zero divisor if there is $b\neq 0$ in $R$ such that $ab=0$.

\begin{thm}
Let $\mathcal{H}$ be a $d$-uniform $r$-partite ($2\leq d\leq r$) hypergraph which satisfies the condition $(\ast\ast)$ for $r\geq 2$. Then $\mathcal{H}$ is unmixed if and only if for every $1\leq q\leq n$, and every $1\leq i_{1}<i_{2}< \ldots <i_{r-d+2}\leq r$, $\overline{x}_{qi_{1}}+\overline{x}_{qi_{2}}+ \ldots +\overline{x}_{qi_{r-d+2}}$ is not a zero divisor in $K[\mathcal{H}]$. Here $\overline{x}_{qi_{t}}$ denotes the image of $x_{qi_{t}}$ in $K[\mathcal{H}]$. 
\begin{proof}
Let $\mathcal{H}$ be unmixed. If for some $1\leq q\leq n$ and some $1\leq i_{1}<i_{2}< \ldots <i_{r-d+2}\leq r$, $\overline{x}_{qi_{1}}+\overline{x}_{qi_{2}}+ \ldots +\overline{x}_{qi_{r-d+2}}$ is zero divisor in $K[\mathcal{H}]$, then there is a polynomial $f\notin I(\mathcal{H})$ in 
\[S=K[x_{11}, \ldots , x_{n1}, x_{12}, \ldots , x_{n2}, \ldots , x_{1r}, \ldots , x_{nr}]\]
such that $f.(x_{qi_{1}}+x_{qi_{2}}+ \ldots +x_{qi_{r-d+2}})\in I(\mathcal{H})$. The ideal $I(\mathcal{H})$ is a monomial ideal and therefore we may assume that $f$ is a monomial and then each monomial of  must belong to $I(\mathcal{H})$ (see \cite{Her 1}). That is, each monomial of the above polynomial must be divided by some generator of $I(\mathcal{H})$ which comes from a hyperedge. Let $fx_{qi_{t}}$ be such a monomial. Then there is a hyperedge $e_{t}$ in $\mathcal{H}$ such that $X_{e_{t}}|fx_{qi_{t}}$. But $X_{e_{t}}\nmid f$. Then $x_{qi_{t}}|X_{e_{t}}$ and $e_{t}\setminus \{x_{qi_{t}}\}$ is a subminimal edge $\mathfrak{e}_{t}$ and $X_{\mathfrak{e}_{t}}|f$. Therefore, for $x_{qi_{1}}, x_{qi_{2}}, \ldots, x_{qi_{r-d+2}}$, There are $r-d+2$ submaximal edges $\mathfrak{e}_{1}, \mathfrak{e}_{2}, \ldots, \mathfrak{e}_{r-d+2}$ such that $X_{\mathfrak{e}_{t}}|f$, for $1\leq t\leq r-d+2$, and    
\[\mathfrak{e}_{1}\thicksim x_{qi_{1}}, \mathfrak{e}_{2}\thicksim x_{qi_{2}}, \ldots , \mathfrak{e}_{r-d+2}\thicksim x_{qi_{r-d+2}}.\]
Now by theorem 3.3, the set
\[\mathfrak{e}_{1}\cup \mathfrak{e}_{2}\cup\ldots \cup \mathfrak{e}_{r-d+2}\]
contains a hyperedge $e$. Now $X_{\mathfrak{e}_{1}\cup \mathfrak{e}_{2}\cup\ldots \cup \mathfrak{e}_{r-d+2}}|f$. Then $X_{e}|f$. Then $f\in I(\mathcal{H})$, a contradiction.

Conversely, let for for every $1\leq q\leq n$, and every $1\leq i_{1}<i_{2}< \ldots <i_{r-d+2}\leq r$, $\overline{x}_{qi_{1}}+\overline{x}_{qi_{2}}+ \ldots +\overline{x}_{qi_{r-d+2}}$ is not zero divisor in $K[\mathcal{H}]$. If $\mathcal{H}$ is not unmixed, by theorem 3.3, there is an integer $1\leq q\leq n$, and submaximal edges $\mathfrak{e}_{1}, \mathfrak{e}_{2}, \ldots, \mathfrak{e}_{r-d+2}$, such that
\[\mathfrak{e}_{1}\thicksim x_{qi_{1}}, \mathfrak{e}_{2}\thicksim x_{qi_{2}}, \ldots , \mathfrak{e}_{r-d+2}\thicksim x_{qi_{r-d+2}}.\]
where where $i_{1}, i_{2} \ldots , i_{r-d+2}$ are distinct and 
\[\mathfrak{e}_{1}\cup \mathfrak{e}_{2}\cup\ldots \cup \mathfrak{e}_{r-d+2}\]
is an independent set. Set $e_{t}=\mathfrak{e}_{t}\cup x_{qi_{t}}$, for $1\leq t\leq r-d+2$. $e_t$'s are hyperedge. Let $X=X_{\mathfrak{e}_{1}\cup \mathfrak{e}_{2}\cup\ldots \cup \mathfrak{e}_{r-d+2}}$. $X$ is not in $I(\mathcal{H})$ but $X.(x_{qi_{1}}+x_{qi_{2}}+ \ldots +x_{qi_{r-d+2}})\in I(\mathcal{H})$, a contradiction.
\end{proof}
\end{thm}


\begin{thebibliography}{20}
\bibitem{Ber} C. Berge, \textit{Hypergraphs}, Elsevier Science Publishers B. V., 1989.
\bibitem{Hag} H. Haghighi, A generalization of Villarreal's result for unmixed tripartite graphs, \textit{Bull. Iranian Math. Soc.}, \textbf{40} (2014), no. 6, 1505-1514.
\bibitem{Her 1} J. Herzog and T. Hibi, \textit{Monomial Ideals}. Springer, 2010.
\bibitem{Her 2} J. Herzog, T. Hibi, H. Ohsugi, Unmixed bipartite graphs and sublattices of Boolian lattices, \textit{J. Algebraic Combin.}, \textbf{30} (2009), no. 4, 415-420.
\bibitem{Jaf} R. Jafarpour-Golzari and R. Zaare-Nahandi, Unmixed r-partire graphs, to appear in Bull. Iranian Math. Soc., arXiv: 1511.00228v1 [math. CO].
\bibitem{Kia} D. Kiani and S. Saeedi Madani, The edge ideal of complete multipartite hypergraphs, arXiv: 1402.2469v1 [math. CO].
\bibitem{Mor} S. Morey, E. Reyes and R. H. Villarreal, Cohen-Macaulay, shellable and unmixed clutters with a perfect matching of K\"{o}nig type, \textit{J. Pure Appl. Algebra}, \textbf{212} (2008), no. 7, 1770-1786.
\bibitem{Rav} G. Ravindra, Well-covered graphs, \textit{J. Combin. Inform.}, System Sci. \textbf{2} (1977), no. 1, 20-21.
\bibitem{Sta} R. Stanly, \textit{Combinatorics and Commutative Algebra}, 2nd Ed., Progress in Math., Birkhauser, 1996.
\bibitem{Vil 1} R. H. Villarreal, \textit{Monomial Algebras}, Marcel Dekker, Inc., New York, 2001.
\bibitem{Vil 2} R. H. Villarreal, Unmixed bipartite graphs, \textit{Rev. Colombiana
Mat.}, \textbf{41} (2007), no. 2, 393-395.
\bibitem{Wes} D. West, \textit{Introduction to graph theory}, 2nd Ed., Prentic Hall, Upper Saddle River, NJ, 2001.
\bibitem{Zar} R. Zaree-Nahandi, Pure simplicial
complexes and well-covered graphs, \textit{Rocky Mountain Journal
of Mathematics}, \textbf{45} (2015), no. 2, 695-702.
\end{thebibliography}
\end{document}